\documentclass[11pt,reqno]{amsart}
\usepackage{amssymb,latexsym}
\usepackage{graphicx}
\usepackage{fancyhdr,amssymb}
\usepackage{url}		

\theoremstyle{plain}
\numberwithin{equation}{section}

\begin{document}
\fancyhead{}
\renewcommand{\headrulewidth}{0pt}
\fancyfoot{}
\fancyfoot[LE,RO]{\medskip \thepage}

\setcounter{page}{1}

\title[Fibonacci Numbers and Identities II]{ Fibonacci Numbers and Identities II}
\author{Cheng Lien Lang}
\address{Department Applied of Mathematics\\
                I-Shou University\\
                Kaohsiung, Taiwan\\
                Republic of China}
\email{cllang@isu.edu.tw}
\thanks{}
\author{Mong Lung Lang}
\email{lang2to46@gmail.com}



\maketitle

\vspace{-.8cm}

\section {Introduction}
\noindent
Let  $a,b, p,q \in \Bbb C$, $q\ne 0$.
Define the generalised Fibonacci sequence $\{W_n \} = \{ W_n(a,b\, ; p,q)\}$ by
 $W_0=a$, $W_1 = b$,
$$
W_n= pW_{n-1} -qW_{n-2}.\eqno(1.1)$$

 \medskip
\noindent
Obviously the definition can be extended to negative subscripts ;  that is,
 for $n = 1, 2, 3,\cdots, $ define
  $$W_{-n} = (pW_{-n+1} -W_{-n+2})/q.\eqno (1.2)$$

 \medskip
 \noindent
 In the case $a=0$, $b=1$, following [Ho], we shall denote the sequence $\{ W_n
  (0,1\,;\, p,q)\}$ by
 $\{u_n\}$. Equivalently,
$$u_n =  W_n(0,1\, ; p,q).\eqno(1.3)$$

\medskip
\noindent
Note that
 $u_{-n} = -q^{-n} u_n$.
 In this article, we will give a
 unified proof of various identities involving $W_n$. The  alternative proof  presented
  in this article (for Melham's, Howard's   and Horadam's identities) uses nothing but
  {\em recurrence} which is slightly different from the original proof.

\section {Recurrence Relation}

\noindent {\bf Lemma 2.1.} {\em  $ A(n) = W_{2n}$, $B(n) =W_{n}W_{n+r}$ and $C(n) = q^n$  satisfy the following
  recurrence relation}
 {\small  $$X(n+3) = (p^2-q)X(n+2) +(q^2-p^2q)X(n+1)+q^3X(n).\eqno (2.1)$$}

   \noindent {\em Proof.}
   One sees easily that $q^n$ satisfies (2.1).
    The following shows that $W_{n}^2 $ satisfies (2.1).  By (1.1),
 {\small   $$W_{n+3}^2 = (pW_{n+2} - qW_{n+1})^2 = p^2W_{n+2}^2 + q^2W_{n+1}^2 - 2pqW_{n+2}W_{n+1},\eqno(2.2)$$

  $$ 2pqW_{n+2}W_{n+1}
  = qW_{n+2}(W_{n+2}+qW_{n})
    +pq(pW_{n+1}-qW_n)W_{n+1}.\eqno(2.3)$$}

 \noindent    Replace the last quantity of the right hand side of (2.2) by (2.3), one has   the following.

   {\small   $$ W_{n+3}^2 =
     (p^2-q)W_{n+2}^2 +(q^2 -p^2q)W_{n+1}^2 - q^2 W_{n+2}W_n +pq^2W_{n+1}W_n.\eqno(2.4)$$}

\noindent By (1.1), {\small  $q^2 W_{n+2}W_n -pq^2W_{n+1}W_n = -q^3W_n$.}  This completes the proof of the
 fact that
 $W_n^2 $ satisfies the recurrence relation (2.1).
 One can show similar to the above that $W_{2n}$ and $ W_{n}W_{n+r}$ satisfy  the recurrence (2.1).
  \qed

 \medskip
   \noindent
  {\bf Lemma 2.2.} {\em Suppose that $x(n)$ and $y(n)$ satisfy the recurrence $(2.1)$. Let $r, s \in \Bbb Z$. Then $ x(n) \pm y(n)$ and  $rx(n+s)$ satisfy the recurrence $(2.1)$.}

\medskip
\noindent  {\bf Lemma 2.3.} {\em Suppose that both $A(n)$ and $B(n)$ satisfy $(2.1). $
  Then $A(n) = B(n)$ if and only if $A(n) = B(n) $ for $n = 0, 1, 2.$}

  \medskip
  \noindent {\em Proof.} Since both $A(n)$ and $B(n)$ satisfy the same recurrence relation,
   $A(n)= B(n)$ if and only if they satisfy the same initial condition. This
    completes the proof of the lemma. \qed

 \medskip
 \noindent {\bf Lemma 2.4.}  {\em Let $W_n$ and $u_n$ be given as in $(1.1)$ and $(1.3)$. Then
  the  following hold.}
  $$  W_r = u_rW_1 -qu_{r-1}W_0,\eqno(2.5)$$
   $$  W_r = u_{r-1}W_2 -qu_{r-2}W_1,\eqno(2.6)$$
     $$  W_r = u_{r-2}W_3 -qu_{r-3}W_2.\eqno(2.7)$$

 \medskip
  \noindent
   {\em Proof.} We note first that (1.1) and (1.2) give the same recurrence relation.
   Denoted by $A(r)$ and $B(r)$ the left and right hand side of (2.5)-(2.7).
    Since  $W_r$, $u_{r}$, $u_{r-1}$  and $u_{r-2}$
   (as functions in $r$) satisfy the recurrence relation
   (1.1), both $A(r)$ and $B(r)$ satisfy (1.1). Hence $A(r) = B(r)$ if and only if they satisfy the same initial
    condition. (2.5)-(2.7) can now be verified with ease. \qed

  \section{On Melham's Identity}

\noindent {\bf Theorem 3.1.} (Melham [M])
{\em   Let $W_n$ be given as in $(1.1)$ and let  $e= pab-qa^2-b^2$. Then
$$W_{n+1}W_{n+2}W_{n+6} -W_{n+3}^3=
eq^{n+1}(p^3W_{n+2} -q^2W_{n+1}).\eqno(3.1)$$}
\noindent  In  the case $p=1, q=-1, a=b=1$, (3.1) becomes
$$F_{n+1}F_{n+2}F_{n+6} -F_{n+3}^3 = (-1)^nF_n,\eqno(3.2)$$
where $F_n$ is the $n$-th Fibonacci number.
In this short note, we give a direct proof of Theorem 3.1
 which uses only the simple  fact that $W_nW_{n+r}$ and $q^n$ satisfy the  recurrence relation (2.1).

  \medskip
  \noindent {\bf Lemma 3.2.} {\em $W_{n+2}W_{n+4} - W_{n+3}^2 = e q^{n+2}$ and
   $W_{n+1}W_{n+6}-W_{n+3}W_{n+4}
   =eq^{n+1}(p^3-pq).$
  }

  \medskip
  \noindent {\em Proof.} Let $A(n) =W_{n+2}W_{n+4} - W_{n+3}^2$, $ B(n) = e q^{n+2}.$
   By Lemmas 2.1 and 2.2,  both $A(n)$ and $B(n)$ satisfy  the recurrence relation (2.1).
    Note that $A(n) = B(n)$ for $n = 0,1,2$. It follows from Lemma 2.3 that
    $A(n) = B(n)$ for all $n$. The second identity can be
     proved similarly.
       \qed

  \medskip
  \noindent {\bf Proof of Theorem 3.1.}
  By the first identity of Lemma 3.2, the left hand side of (3.1) can be written as
   $$
 W_{n+1}W_{n+2}W_{n+6} -W_{n+3}^3=
   W_{n+2}(W_{n+1}W_{n+6}-W_{n+3}W_{n+4})
   +eq^{ n+2}W_{n+3}.\eqno(3.3)$$
 Replace the left hand side of (3.1) by identity (3.3), one has  (3.1) is true if and only if
 $$  W_{n+2}(W_{n+1}W_{n+6}-W_{n+3}W_{n+4})
   +eq^{ n+2}W_{n+3}
   =eq^{n+1}(p^3W_{n+2} -q^2W_{n+1}).\eqno(3.4)
   $$
 By(1.1), $W_{n+3} = pW_{n+2}-qW_{n+1}$ (note that (1.1) and (1.2) give the
  same recurrence relation). Hence   (3.4) is equivalent to
 $$  W_{n+2}(W_{n+1}W_{n+6}-W_{n+3}W_{n+4})
   =eq^{n+1}(p^3-pq)W_{n+2}.\eqno(3.5)
   $$
One may now apply the second identity of Lemma 3.2 to finish our proof of Theorem 3.1. \qed


\section {On Howard's Identity}

\noindent
Let $p$ and $q$ be given as in (1.1).
 Define similarly the sequence $\{V_n\}$ by the same recurrence relation
$$
V_n= pV_{n-1} -qV_{n-2}\,, \,\,V_0 = c, \,\,V_1 = d.\eqno(4.1)$$

\medskip
\noindent A careful study of the proof of Lemma 2.1 reveals that

\begin{enumerate}
\item[(i)] the fact we used to show that
 $W_n^2$ satisfies the recurrence (2.1) is the recurrence (1.1),
 \item[(ii)]  the proof of Lemma 2.1
  is independent of the choice of the values of  $W_0$ and $W_1$.
  \end{enumerate}

  \noindent
   As a consequence,
   the following lemma is clear.

\smallskip
\noindent {\bf Lemma 4.1.}
 {\em
 Let $k , r\in \Bbb Z$. Then
 $rW_nV_{n+k}$, $rW_n u_{n+k}$ and $ rV_{n}u_{n+k}$
  satisfy the recurrence relation $(2.1)$.}

\medskip
\noindent {\bf Theorem 4.2.} (Howard [Ho])
{\em Let $W_n, V_n$ and $ u_n$ be given as in $(1.1)$, $(1.3)$ and $(4.1)$. Then }

 $$W_n^2 -q^{n-j}W_j^2
 = u_{n-j}(bW_{n+j}-qaW_{n+j-1}),\eqno(4.2)$$
 $$W_{m+n+1} = W_{m+1}u_{n+1}-qW_mu_n,\eqno(4.3)$$
 $$V_{m+k}W_{n+k}-q^kV_mW_n =
  u_k(bV_{m+n+k}-qaV_{m+n+k-1}),\eqno(4.4)$$
   $$ W_{n+k}^2 -q^{2k}W_{n-k}^2
   = u_{2k}(bW_{2n} -qaW_{2n-1}).\eqno(4.5)$$

 \medskip
 \noindent  In [Ho], proof of Theorem 4.2 involves mathematical induction, a
   result of Bruckman [B] and a  very neat combinatorial argument.
   In this article, we give  an alternative proof which is direct and uses Lemmas 2.1-2.4 and  4.1 only.
  \medskip

  \noindent
   {\bf  Proof of Theorem 4.2.}  We shall prove (4.3). The rest can be proved similarly.
   Note first that (4.3) can be proved by induction as well.
   Let $m+n+1 = r$.  Since $u_{-n} = -q^{-n} u_n$
    and $q\ne 0$,
    (4.3) becomes
   $$
  -q^{m-r} W_r = W_{m+1}u_{m-r} -W_m u_{m+1-r}.\eqno(4.6)$$

  \medskip
  \noindent Let $A(m) = - q^{m-r}W_r$ and  $ B(m) =W_{m+1}u_{m-r} -W_m u_{m+1-r}.$
   View both $A(m)$ and $B(m)$ as functions in $m$. By Lemmas 2.1 2.2 and  4.1, they both satisfy the
    recurrence relation (2.1). By Lemma 2.3,  $A(m) = B(m)$ if and only if $A(m)= B(m)$
     for $m = 0, 1$ and 2. Equivalently, (4.3) is true if and only if
     $A(0) = B(0)$, $A(1) = B(1)$ and $A(2) = B(2)$. Equivalently,

    $$-q^{-r} W_r= W_1u_{-r} - W_0u_{1-r},   \eqno(4.7)  $$
  $$  - q^{1-r}W_r = W_2u_{1-r} - W_1u_{2-r},  \eqno(4.8)$$
  $$  -q^{2-r}W_r = W_3u_{2-r} - W_2u_{3-r}.
    \eqno(4.9)
    $$

 \medskip
 \noindent
  Since $u_{-n} = -q^{-n} u_n$ and $q \ne 0$,  (4.7), (4.8) and (4.9) of the above  can be written as the following.
  $$  W_r = u_rW_1 -qu_{r-1}W_0,\eqno(4.10)$$
  $$  W_r = u_{r-1}W_2 -qu_{r-2}W_1,\eqno(4.11)$$
    $$  W_r = u_{r-2}W_3 -qu_{r-3}W_2.\eqno(4.12)$$

 \medskip
  \noindent Hence (4.3) is true if and only if (4.10)-(4.12) are true. Applying Lemma 2.4,
   we conclude that (4.3) is true. \qed

\section{Horadam's Identities}
  Horadam [H] studied the fundamental arithmetical properties of  $\{W_n\}$
and provided us with many interesting identities. A careful checking of his identities suggested that
 Identities (3.14), (3.15), (4.1)-(4.12) and (4.17)-(4.22) of [H] can be
  verified by our method as the functions involved in those identities satisfy the recurrence (2.1)
   of the present paper.

  \section{Discussion}
  The technique we presented in the proof of Theorems 3.1 and  4.2 can be
   applied to the verification of identities of the form $A(n)= B(n)$ if one can show

  \begin{enumerate}
   \item[(i)]   $A(n)$ and $B(n)$ satisfy the same
    recurrence relation,
    \item[(ii)] $A(n)$ and $B(n)$ satisfy the same initial condition.
\end{enumerate}

\noindent The verification of (i) is not difficult (see Appendix B of [LL]). (ii) can be checked with ease.
Note that it is not just our purpose to reprove the identities but to illustrate the importance and usefulness
 of the
 {\em recurrence} (such as (2.1) of the present paper)

\medskip

\noindent MSC2010: 11B39, 11B83
\medskip


\end{document}